# On the Flexibility of Suspensions

Gerald D. Nelson
nelso229@comcast.net

**Abstract**

We study the flexibility of suspensions (polyhedra having the combinatorial structure of dipyramids) that have an even number of vertexes and provide arguments that there are least five distinct types of flexible suspensions.

## 1. Introduction

Suspensions are polyhedra having the combinatorial structure of dipyramids that are comprised of two polyhedral caps each having N triangular faces as illustrated in Fig. 1. The two apical vertexes **u** and **w**, are of index N while the non-apical vertexes at the base of the caps, $v_1, v_2, \ldots v_N$, are all of index 4.

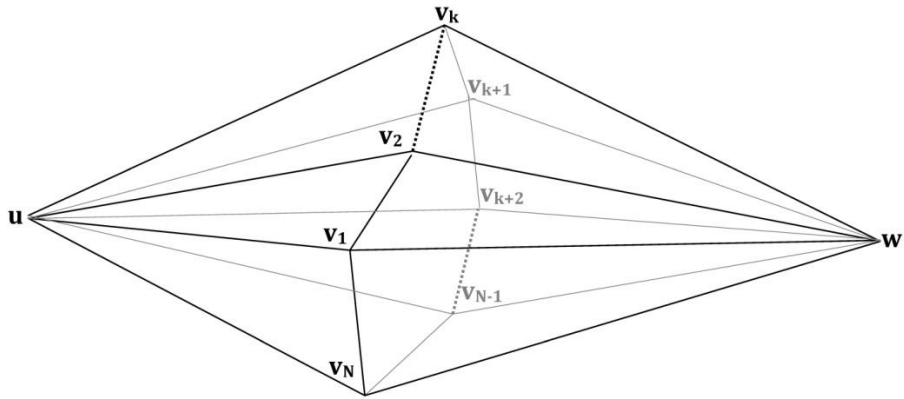

**Figure 1.** Suspension with apical vertexes of Index N

The flexibility of these polyhedra has been established for N=4, Bricard octahedra [1], N=6 [2] and for specific constructions for even N=6..32 [3]. Flexible suspensions for odd N are not known. Additionally this topic has been widely studied; see [4], [5] and [6] for example.

In this paper we offer arguments that there are at least five types of flexible suspensions for even N>4. We proceed as follows. Sec. 2 describes notation which we use to characterize suspensions. In Sec. 3 we address the flexibility of five specific types of suspensions having an even number of vertexes. Sec. 4 contains theorems regarding the volume of these suspensions and the "strength" of their flexibility. Sec. 5 contains concluding remarks and identifies some flexible suspensions that are not addressed in this paper.

## 2. Notation

This section introduces notation, taken from [3], that we use to characterize edge lengths, face angles and dihedral angles that appear in suspensions. With reference to Fig. 1: N is an even integer that always refers to the number of non-apical vertexes and N=2M for some M>2. The notation used for edge lengths, face angles an dihedral angles is shown in Table 1. For all items k=1..N, and here and elsewhere, the cyclic convention $(N+1) \equiv 1$ holds.

| Edge lengths | | Face angles | | Dihedral angles | |
|---|---|---|---|---|---|
| $l_k$ | $|u-v_k|$ | $\alpha_k$ | $\angle v_k u v_{k+1}$ | $\varepsilon_k$ | $\angle v_k v_{k+1}$ |
| $m_k$ | $|w-v_k|$ | $\beta_k$ | $\angle v_{k+1} v_k u$ | $\Delta_k$ | $\angle w v_k$ |
| $L_k$ | $|v_k-v_{k+1}|$ | $\gamma_k$ | $\angle u v_{k+1} v_k$ | $\delta_k$ | $\angle u v_k$ |
| | | $A_k$ | $\angle v_{k+1} w v_k$ | | |
| | | $B_k$ | $\angle w v_k v_{k+1}$ | | |
| | | $\Gamma_k$ | $\angle v_k v_{k+1} w$ | | |

**Table 1**. Notation Definition

Fig. 2 illustrates this notation for two typical faces, $uv_{k+1}v_k$ and $wv_kv_{k+1}$. By convention this is the view of the outsides of the faces. The vectors $(v_{k+1}-u) \times (v_k-u)$ and $(v_k-w) \times (v_{k+1}-w)$, where **x** is the vector cross product operator, determine the direction vectors (pointing into the page) associated with these faces.

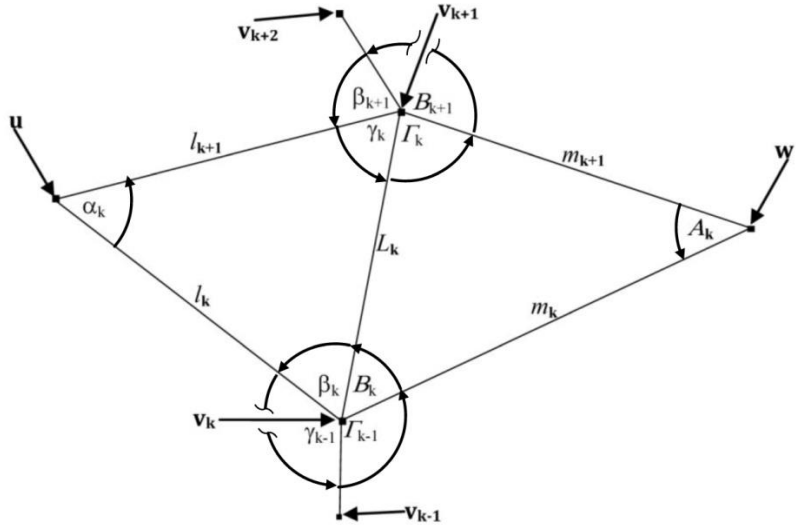

**Figure 2.** Notation.

## 3. Flexible Suspensions

In this section we provide arguments that suspensions with an even number of vertexes and that are identified in [3] by five generalizations of Bricard octahedra identified as types I-OEE, II-AEE, II-OEE, III-OAE and III-OAS are flexible. In these type identifiers the Roman numeral refers to the type of the related Bricard octahedron; the other letters are acronyms for geometric characteristics: OEE – Opposite Edges Equal, AEE – Adjacent Edges Equal, OAE – Opposite Angles Equal and OAS – Opposite Angles Supplementary.

Preparatory to a discussion of individual suspension types we define a variable of flexion and a coordinate model that can be applied to all five types and is useful in proofs of the flexibility of the specific types. The variable of flexion is taken to be the distance between the apical vertexes: z = |**u** − **w**|. The coordinate model is recursive with respect to the manner in which non-apical vertexes are defined; see Table 2.

| vertex | x-coordinate | y-coordinate | z-coordinate |
|--------|--------------|--------------|--------------|
| **u** | 0 | 0 | $\dfrac{z}{2}$ |
| **w** | 0 | 0 | $\dfrac{-z}{2}$ |
| **v_k** | $r_k \cos\theta_k$ | $r_k \sin\theta_k$ | $\dfrac{m_k^2 - l_k^2}{2z}$ |

**Table 2.** Coordinate Model Vertex Definition

Here $r_k$ and $\theta_k$ are functions of the variable of flexion z and the edge lengths as defined in Sec. 2. $r_k$ is defined to be the closest distance from the vertex **v_k** to the line connecting vertexes **u** and **w**; see Fig. 3:

$$r_k = \sqrt{\frac{(m_k^2 + l_k^2)}{2} - z_k^2 - \frac{z^2}{4}}$$

for k=1..N. The angle $\theta_k$ is the rotation angle about the z-axis with respect to the x-z plane; it is defined recursively:

$$\begin{aligned} \theta_{k+1} &= \theta_k \pm \Delta\theta_k \text{ where} \\ \Delta\theta_k &= \cos^{-1} t_k \text{ and} \\ t_k &= \frac{r_{k+1}^2 + r_k^2 - L_k^2 + (z_{k+1} - z_k)^2}{2 r_k r_{k+1}} \end{aligned} \quad (1)$$

for k=1..N-1, for an appropriately defined $\theta_1$ and choice of the signs associated with $\Delta\theta_k$.

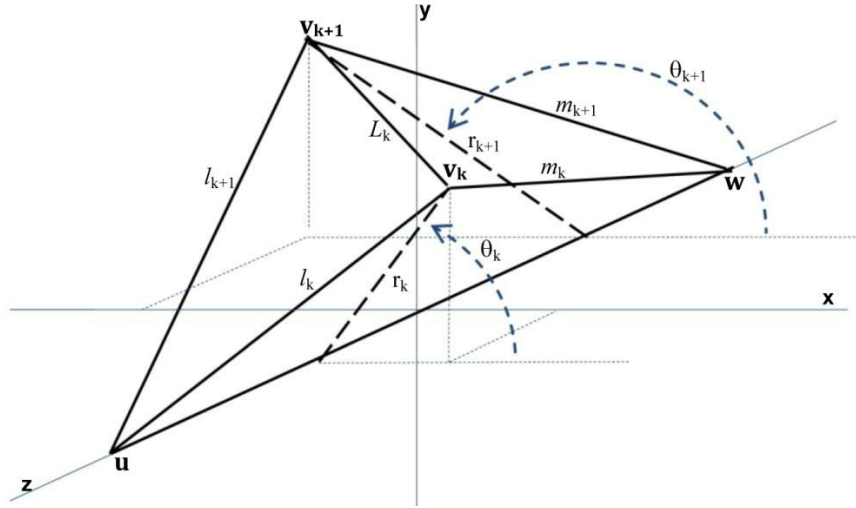

**Figure 3.** Coordinate Model.

The use of this coordinate model is advantageous since it can be shown by direct evaluation that the following theorem holds:

**Theorem I.** Suspensions defined by Table 2 are flexible if the edge length given by $|v_N - v_1|$ is not a function of variable of flexion z.

Proof: By direct evaluation it is seen that the following is true

$$|v_k - u| = l_k$$
$$|v_k - w| = m_k \quad \text{for k=1..N and}$$
$$|v_{k+1} - v_k| = L_k \text{ for k=1..N-1}.$$

Since these results are independent of the variable z it is apparent that any choice of edge lengths, $\theta_1$ and signs associated with $\Delta\theta_k$ for which the equation

$$|v_N - v_1| = L_N \tag{2}$$

is true and is independent of the variable z effectively defines a flexible suspension. ∎

Additionally it is easy to see that the variable of flexion z is related to the each of the dihedral angles $\varepsilon_k$ by the relationship for k=1..N:

$$z^2 = l_k^2 + m_k^2 - 2m_k l_k (\cos\beta_k \cos B_k + \sin\beta_k \sin B_k \cos\varepsilon_k). \tag{3}$$

We utilize these results in two ways: for the types I-OEE, II-AEE and II-OEE a complete set of edge lengths is introduced into the coordinate model, Table 2, and Theorem I is shown to be satisfied for specific choices of $\theta_1$ and signs in equations (1). For types III-OAE and III-OAS a more constructive approach is used as these types are not characterized by a complete set of edge lengths. A recursive definition of edge lengths and associated face angles that leads to flexible suspensions is developed.

The type **I-OEE** suspension is a generalization of Bricard octahedra of the first type and is characterized by 3M independent edge length parameters which satisfy the relationships:

$$m_k = l_{k+M},$$
$$m_{k+M} = l_k \text{ and} \tag{4}$$
$$L_{k+M} = L_k \text{ for k=1..M}.$$

For this type we assert:

**Theorem II.** Type I-OEE suspensions are flexible.

Proof. For the parameters defined in equations (4) it is easy to see that when $\theta_1$ is defined by:

$$\theta_1 = \frac{\pi - \Sigma_M}{2}$$

where

$$\Sigma_M = \sum_{k=1..M} \Delta\theta_k$$

we have $\theta_{k+M} = \pi - \theta_k$ and the coordinate model in Table 2 simplifies to:

$$x_{k+M} = -x_k,$$
$$y_{k+M} = y_k \text{ and}$$
$$z_{k+M} = -z_k \text{ for k=1..M.}$$

By direct evaluation it is seen that equation (2) is satisfied and is independent of the variable z. ∎

Before moving to the next suspension type it is noted that the above coordinate model simplification illustrates the axial symmetry that exists for this type and which carries over from the axial symmetry of Bricard octahedra of the first type. Here we have two apical vertexes that move at right angles with respect to the axis of symmetry and N non-apical vertexes that pair-wise move symmetrically about that axis, in this case the y-axis.

The type **II-AEE** suspension is a generalization of Bricard octahedra of the second type and is characterized by 3M independent edge length parameters that satisfy the relationships:

$$m_1 = l_1,$$
$$m_k = l_{N-k+2} \text{ for k=2..N and} \tag{5}$$
$$L_{k+M} = L_{M-k+1} \text{ for k=1..M.}$$

For this type we assert:

**Theorem III.** Type II-AEE suspensions are flexible.

Proof. For the parameters defined in equations (5) it is easy to see that when $\theta_1 = 0$ the coordinate model in Table 2 simplifies to:

$$x_1 = r_1,$$
$$y_1 = 0,$$
$$z_1 = 0,$$
$$x_{k+M} = x_{M-k+2},$$
$$y_{k+M} = y_{M-k+2} \text{ and}$$
$$z_{k+M} = -z_{M-k+2} \text{ for k=2..M.}$$

For k=M it is seen that

$$x_N = x_2,$$
$$y_N = y_2 \text{ and}$$
$$z_N = -z_2$$

thus by direct evaluation, equation (2) is satisfied and is independent of the variable z. ∎

It is noted that the above coordinate model simplification illustrates a form of planar symmetry that exists for this type and which carries over from the planar symmetry of Bricard octahedra of the second type. The apical vertexes and all other vertexes move symmetrically with respect to the x-y plane except $\mathbf{v_1}$ and $\mathbf{v_{M+1}}$ which lie in the plane.

The type **II-OEE** suspension is a second generalization of Bricard octahedra of the second type and is characterized by 3M independent edge length parameters and that satisfy the relationships:

$$l_{k+M} = l_k,$$
$$m_{k+M} = m_k \text{ and} \tag{6}$$
$$L_{k+M} = L_k \text{ for k=1..M.}$$

For this suspension type we assert:

**Theorem IV.** Type II-OEE suspensions are flexible.

Proof. For the parameters defined in equations (6) it is easy to see that when $\theta_1$ is defined by

$$\theta_1 = \frac{\pi - \Sigma_M}{2}$$

where

$$\Sigma_M = \sum_{k=1..M} \Delta\theta_k$$

the coordinate model in Table 2 simplifies to:

$$x_{k+M} = -x_k,$$
$$y_{k+M} = y_k \text{ and}$$
$$z_{k+M} = z_k \text{ for } k=1..M.$$

By direct evaluation it is seen that equation (2) is satisfied and is independent of the variable z. ∎

The above coordinate model simplification illustrates a second form of planar symmetry that exists for this type and which carries over from the planar symmetry of Bricard octahedra of the second type.

The suspension type **III-OAE** is the first of two generalizations of Bricard octahedra of the third type and is characterized by angular relationships that permit two positions in which all vertexes are co-planar. The first of these relationships is an open folding in which $\delta_1=\delta_L=0$ and $\delta_k=\pi$ for k=2..(L-1) and k=(L+1)..N for some L where 2<L<(N-2). The second is a compact folding in which $\delta_1=\delta_L=\pi$ and $\delta_k=0$ for k=2..(L-1) and k=(L+1)..N. It is easy to see that for these foldings the following relationships hold for the face angles at the apical vertex **u**:

$$\sum_{k=1..k_x} a_{2k-1} = \sum_{k=k_x+1..M} a_{2k-1}$$

and when L is even ($L = 2k_x$) (7)

$$\sum_{k=1..k_x-1} a_{2k} = \sum_{k=k_x..M} a_{2k}$$

or when L is odd ($L = 2k_x+1$)

$$\sum_{k=1..k_x} a_{2k} = \sum_{k=k_x+1..M} a_{2k}$$

Face angles $A_k$ at the apical vertex **w** satisfy identical equations. Additionally vertexes **v₁** and **v_L** are OAS while the remaining N-2 non-apical vertexes are OAE.

This type is parameterized by the parameter set $\{l_1,m_1,l_2,m_2,L_1,L_3,.. L_{2k-1},…, L_{N-3}\}$ which effectively defines the first two triangular faces of the suspension along with odd indexed edge lengths at the non-apical vertexes while the remaining parameters are given by recursive definitions of edge lengths and face angles where each stage of the recursion defines four faces of the suspension in terms of the definitions from the previous stage, or in the case of the first stage from the parameterization given above.

At the k-th stage, for k=1..(M-1), the edges $l_{2k}$ and $m_{2k}$ and angles $\beta_{2k}$ and $B_{2k}$ are available from previous stages; the edge length $L_{2k+1}$ is a parameter for k<(M-1) while for the last stage (k=M-1) the apical angles $\alpha_{N-1}$ and $\alpha_N$ are defined by equations (7).

Additional parameters of the four faces at the vertex $\mathbf{v_{2k+1}}$ ($\mathbf{uv_{2k+1}v_{2k}}$, $\mathbf{wv_{2k}v_{2k+1}}$, $\mathbf{uv_{2k+2}v_{2k+1}}$, $\mathbf{wv_{2k+1}v_{2k+2}}$) are defined by law of sines and other relationships as described in the following discussion.

The angle $\beta_{2k+1}$ is determined from solutions of the quadratic equation:

$$a_k X(\beta_{2k+1})^2 + b_k X(\beta_{2k+1}) + c_k = 0$$

where the function $X(\phi)$ is defined by

$$X(\phi) = \cotan\left(\frac{\phi}{2}\right)$$

The coefficients $a_k$, $b_k$ and $c_k$ and the angle $B_{2k+1}$ have two definitions depending upon two cases of interest; either

$$a_k = R_k(B_k + C_k R_k),$$
$$b_k = -2A_k R_k \text{ and}$$
$$c_k = -(C_k + B_k R_k) \text{ with}$$
$$R_k = \frac{(1 - K_k)}{(1 + K_k)} \text{ and}$$
$$X(B_{2k+1}) = R_k X(\beta_{2k+1}).$$

or

$$a_k = B_k R_k - C_k$$
$$b_k = -2A_k R_k \text{ and}$$
$$c_k = -R_k(B_k - C_k R_k) \text{ with}$$
$$R_k = \frac{(1 + K_k)}{(K_k - 1)} \text{ and}$$
$$X(B_{2k+1}) = \frac{R_k}{X(\beta_{2k+1})}$$

where $A_k$, $B_k$ and $C_k$ are defined by:

$$A_k = l_{2k} \cos\beta_{2k} - m_{2k} \cos B_{2k}$$
$$B_k = m_{2k} \sin B_{2k}$$
$$C_k = -l_{2k} \sin\beta_{2k}$$

In both cases $K_k$ is derived from a relationship that exists between each of the following vertex pairs:

$$(\mathbf{v_1},\mathbf{v_3}),(\mathbf{v_2},\mathbf{v_5}),(\mathbf{v_4},\mathbf{v_7}),\ldots (\mathbf{v_{2k}},\mathbf{v_{2k+3}}),\ldots (\mathbf{v_{N-4}},\mathbf{v_{N-1}}), (\mathbf{v_{N-2}},\mathbf{v_N}) \text{ for k=1..M-2.}$$

These relationships are a generalization of those observed in Bricard octahedra of the third type, in the flexible suspension defined in [2] and in examples of III-OAE flexible suspensions [3]. They are defined by equations involving the dihedral angles at the vertex pairs. In terms of the functions:

$$V_P(\delta, \varepsilon) = \tan\left(\frac{\delta}{2}\right) \tan\left(\frac{\varepsilon}{2}\right)$$

and

$$V_R(\delta, \varepsilon) = \frac{\tan\left(\frac{\delta}{2}\right)}{\tan\left(\frac{\varepsilon}{2}\right)}$$

for each of the pairs ($v_i$, $v_j$) it is observed that for adjacent dihedral angles that the values are equal in magnitude (but not necessarily of the same sign). For our purposes here for any pair of vertex indexes, say (i,j), as defined above then for III-OAE vertexes

$$V_P(\delta_i, \varepsilon_i) = \pm V_P(\delta_j, \varepsilon_j) \tag{8}$$

while the ratio function $V_R(\delta,\varepsilon)$ is used for III-OAS vertexes. For the specific pair of indexes, i=2k-2 and j=2k+1 the value of $K_k$ is $=V_P(\delta_{2k-2},\varepsilon_{2k-2})$ if $v_{2k-2}$ is an III-OAE vertex, otherwise $=V_R(\delta_{2k-2},\varepsilon_{2k-2})$. For the first stage the indexes are i=1 and j=3. $V_P(\delta_{2k-2},\varepsilon_{2k-2})$ and $V_R(\delta_{2k-2},\varepsilon_{2k-2})$ are defined at a value of the variable z from equation (3) for k=1 and $\varepsilon_1 = \frac{\pi}{2}$ or $= \frac{3\pi}{2}$.

With the angles $\beta_{2k+1}$ and $B_{2k+1}$ defined the parameterization is completed as follows; in the case that the vertex $v_{2k+1}$ is OAE:

$$\gamma_{2k} = B_{2k+1} \text{ and}$$
$$\Gamma_{2k} = \beta_{2k+1}$$

otherwise:

$$\gamma_{2k} = \pi - B_{2k+1} \text{ and}$$
$$\Gamma_{2k} = \pi - \beta_{2k+1}$$

In either case:

$$\alpha_{2k} = \pi - \beta_{2k} - \gamma_{2k} \text{ and}$$
$$A_{2k} = \pi - B_{2k} - \Gamma_{2k}$$

The remaining edge lengths and angles are defined by law of sines and law of cosines:

$$\begin{aligned}
L_{2k} &= \frac{l_{2k} \sin \alpha_{2k}}{\sin \gamma_{2k}}, \\
l_{2k+1} &= \frac{l_{2k} \sin \beta_{2k}}{\sin \gamma_{2k}}, \\
m_{2k+1} &= \frac{m_{2k} \sin B_{2k}}{\sin \Gamma_{2k}}, \\
l_{2k+2} &= \sqrt{l_{2k+1}^2 + L_{2k+1}^2 - l_{2k+1} L_{2k+1} \cos \beta_{2k+1}}, \\
m_{2k+2} &= \sqrt{m_{2k+1}^2 + L_{2k+1}^2 - m_{2k+1} L_{2k+1} \cos B_{2k+1}}, \\
\sin \beta_{2k+1} &= \frac{l_{2k+2} \sin \gamma_{2k+1}}{l_{2k+1}}, \\
\sin B_{2k+1} &= \frac{m_{2k+2} \sin \Gamma_{2k+1}}{m_{2k+1}}, \\
\alpha_{2k+1} &= \pi - \beta_{2k+1} - \gamma_{2k+1} \text{ and} \\
A_{2k+1} &= \pi - B_{2k+1} - \Gamma_{2k+1}
\end{aligned} \tag{9}$$

provided that triangle inequalities are satisfied.

Also the signs of the rotation angles in equations (1) for vertexes $v_{2k}$ and $v_{2k+1}$ are determined as a part of the k-th stage. Valid signs yield a geometry in which the above dihedral relationships, equation (8), are maintained. At the last stage additional caveats are imposed upon the recursion; for example, the apical angle $\alpha_N$ derived from equations (7) must $= \pi - \beta_N - \gamma_N$ ($\gamma_N = B_1$); a choice of signs in equations (1) must exist such that the edge length $\lambda_N(z) = |v_N - v_1|$ is the same when derived from either apical angle $\alpha_N$ or apical angle $A_N$; and equation (8) holds the vertex pair ($v_{N-2}$, $v_N$).

From equations (9) it is seen that the k-th stage of the recursive definition effectively defines all parameters associated with the vertexes $v_{2k}$ and $v_{2k+1}$. For the stage with k=1 the first five parameters provide the basis for the definition while for k>1 the k-1 stage definition is sufficient. At each stage a valid solution is determined and the next stage of definition is undertaken. Valid solutions at each stage result in a suspension about which we assert the following theorem:

**Theorem V.** Type III-OAE suspensions are flexible.

Proof. From the definition of this type (equations (7),(8) and (9) and supporting text) it is seen that we have defined a polyhedron which has two non-apical vertexes ($v_1$ and $v_L$) that are OAS while the remainder are OAE.

It is well known that the angular relationships (dihedral and face angles) of vertexes of index 4 are governed by two types of equations. The first of these provides a relationship between adjacent dihedral angles and related face angles; the "equation of tetrahedral angle" [1; Eq. 1]. The second defines the relationship that exists between opposing dihedral angles and can be found in [1] but is not identified there by number.

The first case has solutions [1; Eqs. 4 and 5] which when formulated for dihedral angles $\delta_k$ and $\varepsilon_k$ and when vertex $v_k$ is III-OAE:

when vertex $v_k$ is III-OAE:
$$V_P(\delta_k, \varepsilon_k) = C_R(B_k, \beta_k) \text{ or } = -S_R(B_k, \beta_k),$$
$$V_R(\delta_k, \varepsilon_k) = -C_R(B_k, \beta_k) \text{ or } = S_R(B_k, \beta_k).$$
(10)

Functions $S_R(\rho,\sigma)$ and $C_R(\rho,\sigma)$ are defined by:

$$S_R(\rho, \sigma) = \frac{\sin\left(\frac{\rho - \sigma}{2}\right)}{\sin\left(\frac{\rho + \sigma}{2}\right)} \text{ and}$$

$$C_R(\rho, \sigma) = \frac{\cos\left(\frac{\rho - \sigma}{2}\right)}{\cos\left(\frac{\rho + \sigma}{2}\right)}.$$

The second case when formulated for the non-apical vertex $v_k$ and dihedral angles $\varepsilon_{k-1}$ and $\varepsilon_k$ for k=1..N has the form:

$$\cos\gamma_{k-1}\cos\Gamma_{k-1} + \sin\gamma_{k-1}\sin\Gamma_{k-1}\cos\varepsilon_{k-1} = \cos\beta_k\cos B_k + \sin\beta_k\sin B_k\cos\varepsilon_k$$

while for dihedral angles $\delta_k$ and $\Delta_k$:

$$\cos\gamma_{k-1}\cos\beta_k + \sin\gamma_{k-1}\sin\beta_k\cos\delta_k = \cos\Gamma_{k-1}\cos B_k + \sin\Gamma_{k-1}\sin B_k\cos\Delta_k.$$

For both OAE and OAS vertexes these equations reduce to:

$$\cos\varepsilon_{k-1} = \cos\varepsilon_k \text{ and}$$
$$\cos\delta_k = \cos\Delta_k.$$
(11)

Since neither equations (10) or (11) constrain the dihedral angles to a specific value of the variable z we conclude that the defined suspension is indeed flexible. ∎

The suspension type **III-OAS** is the second generalization of Bricard octahedra of the third type and is also characterized by angular relationships that permit two positions in which all vertexes are co-planar. It is similar to the III-OAE type with three differences. The open folding is circular rather than fan like with $\delta_k=\pi$ for all k. The compact folding has $\delta_k=0$ for all k. All non-apical vertexes are OAE.

**Theorem VI.** Type III-OAS suspensions are flexible.

Proof. The proof is very similar to the previous proof and is left as an exercise. ∎

## 4. Strong Flexibility and Volume

There is a result regarding the nature of the flexibility of suspensions that is worth memtioning. By defining "strong flexibilty" as the property of a flexible polyhedron for which all dihedral angles are non-constant under flexion we can state:

**Theorem IX.** Flexible suspensions of the types I-OEE, II-AEE, II-OEE, III-OAE and III-OAS exhibit strong flexibility.

Proof: From equation (3) it is seen that the variable of flexion $z$ is related to the each of the dihedral angles $\varepsilon_k$ by a variable relationship; clearly the latter are non-constant. Further, the dihedral angles $\delta_k$ and $\Delta_k$ are each related to $\varepsilon_k$ by equations of the form of the "equation of the tetrahedral angle" [1; Eq. 1]; therefore cannot be constant. ∎

Additionally we can state:

**Theorem X.** Flexible suspensions of the types I-OEE, II-AEE, II-OEE, III-OAE and III-OAS have zero volume.

Proof: It can be shown by direct evaluation that each of the types I-OEE, II-AEE and II-OEE have pairs of faces whose contribution to volume exactly cancel one another. For example for type I-OEE the two faces of interest are **uv**$_{k+1}$**v**$_k$ and **wv**$_{k+M}$**v**$_{k+M+1}$.

On the other hand the types III-OAE and III-OAS have two positions in which all faces lie in the same plane, thus have zero volume in these positions; it is well known [7] that volume is a constant in a flexible polyhedron. ∎

## 5. Conclusion

In this paper we have shown that there are five distinct types of flexible suspensions that have an even number of vertexes. It is not known if there are other such basic types, however it is worth mentioning that there are two other types of flexible polyhedra that can be defined based upon the five suspensions we have described. While a full description of these is beyond the scope of this paper the following serves as a brief summary of their features.

The first of these types are compound suspensions that are formed by joining together two of the suspensions of the above five types. The union is accomplished by uniting the individual polyhedra at identical adjacent faces that include both apical vertexes; eg. faces **uv**$_{k+1}$**v**$_k$ and **wv**$_k$**v**$_{k+1}$ for some limited number of contiguous indexes k=k$_1$..k$_2$.

The second of these types are "extended" suspensions that are formed by joining together portions of two or more of the same suspension some of which may be scaled differently thereby yielding flexible polyhedra of either genus 0 or genus 1. The resulting flexible polyhedra, in the genus 0 case, have two apical caps on index N, as do suspensions, that are not joined by their adjacent faces but are separated by angular rings of kN quadrangular faces for some integer k>0. These are identical in form to the "extended" Bricard octahedra described in [8].

As a final note in closing: we claim no success in the determination of flexible suspensions that have an odd number of vertexes.

## Acknowledgements

The author wishes to thank his wife Verla Nelson for her encouragement and support.

**Gerald Nelson** is a retired software engineer; was employed by Honeywell, Inc. (when it was a Minnesota based company) and by MTS Systems Corp. of Eden Prairie, Minnesota. He has a Masters degree in Mathematics from the University of Minnesota.